\documentclass[12pt,leqno]{amsart}

\usepackage{times}
\usepackage[T1]{fontenc}
\usepackage{mathrsfs}
\usepackage{latexsym}
\usepackage[dvips]{graphics}
\usepackage{epsfig}
\usepackage{mathtools}
\usepackage{amsmath,amsfonts,amsthm,amssymb,amscd}
\usepackage{color}
\usepackage{verbatim}
\usepackage{hyperref}
\usepackage{graphicx}
\usepackage{pstricks}
\usepackage{pgf,tikz}
\usepackage{dsfont}
\usetikzlibrary{arrows}
\usepackage{subcaption}
\usepackage{tikz}
\usepackage{pgfplots}
\usepackage{tikz-3dplot}
\tdplotsetmaincoords{60}{115}
\pgfplotsset{compat=newest}
\usepackage{subcaption}
\usetikzlibrary{decorations.pathreplacing}
\usepackage{tkz-euclide}
\newcommand\be{\begin{equation}}
\newcommand\ee{\end{equation}}
\newcommand\benn{\begin{equation*}}
\newcommand\eenn{\end{equation*}}
\newcommand\bc{\begin{center}}
\newcommand\ec{\end{center}}
\newcommand\bea{\begin{eqnarray}}
\newcommand\eea{\end{eqnarray}}
\newcommand\beann{\begin{eqnarray*}}
\newcommand\eeann{\end{eqnarray*}}
\newcommand\bi{\begin{itemize}}
\newcommand\ei{\end{itemize}}
\newcommand\ben{\begin{enumerate}}
\newcommand\een{\end{enumerate}}
\newcommand\bml{\begin{multline}}
\newcommand\eml{\end{multline}}
\newcommand\bmln{\begin{multline*}}
\newcommand\emln{\end{multline*}}

\newtheorem{thm}{Theorem}[section]

\newtheorem{cor}[thm]{Corollary}
\newtheorem{lem}[thm]{Lemma}
\newtheorem{prop}[thm]{Proposition}

\newtheorem{rek}[thm]{Remark}

\newcommand{\h}[1]{\widehat{#1}}
\newcommand{\lla}[1]{\left \{  #1  \right \}}
\newcommand{\pa}[1]{\left (  #1 \right )}
\newcommand{\abs}[1]{\left |  #1  \right |}

\newcommand{\corch}[1]{\left [  #1  \right ] }

\newcommand{\R}{\ensuremath{\mathbb{R}}}

\newcommand{\dps}{\displaystyle}

\newcommand{\df}{\mathrm{d}}
\newcommand{\ve}{\varepsilon}

\newcommand{\tri}{\bigtriangleup}

\begin{document}

\definecolor{qqqqff}{rgb}{0.,0.,1.}
\definecolor{ffqqqq}{rgb}{1.,0.,0.}

\title{A Mattila-Sj\"{o}lin theorem for triangles}

 \author{Eyvindur Ari Palsson}
 \thanks{The work of the first listed author was supported in part by Simons Foundation Grant no. 360560, and the work of the second listed author was supported in part by NSF grant DMS1907435}
 \email{\textcolor{blue}{\href{mailto:eap2@williams.edu} {palsson@vt.edu}}}
 \address{Department of Mathematics, Virginia Tech, Blacksburg, VA 24061}

 \author{Francisco Romero Acosta}
 \email{\textcolor{blue}{\href{mailto:eap2@williams.edu} {jfromero@vt.edu}}}
 \address{Department of Mathematics, Virginia Tech, Blacksburg, VA 24061}




\begin{abstract}
We show for a compact set $E \subset \mathbb{R}^d$, $d \geq 4$, that if the Hausdorff dimension of $E$ is larger than $\frac{2}{3}d+1$, then the set of congruence classes of triangles formed by triples of points of $E$ has nonempty interior. Here we understand the set of congruence classes of triangles formed by triples of points of $E$ as the set
$$\Delta_{\text{tri}}(E) = \left \{ (t,r, \alpha) : |x-z|=t, |y-z|=r \, \text{ and  }\, \alpha= \alpha(x,z,y), \  x,y,z \in E \right \},$$
where $\alpha (x,z,y)$ denotes the angle formed by $x$, $y$ and $z$ , centered at $z$. This extends the Mattila-Sj\"{o}lin theorem that establishes a non-empty interior for the distance set instead of the set of congruence classes of triangles. These theorems can be thought of as refinements and extensions of the statements in the well known Falconer distance problem.
\end{abstract}

\maketitle



\section{Introduction}

The celebrated Falconer's distance problem asks how large does the Hausdorff dimension of a compact set $E \subset \R^d$, $d \geq 2$, need to be to ensure that its distance set $\Delta (E) = \left \{ \left | x-y \right |: x,y \in E \right \}$ has positive Lebesgue measure. Falconer conjectured that if the Hausdorff dimension of $E$ is larger than $\frac{d}{2}$, then the distance set $\Delta (E)$ has positive Lebesgue measure. In 1985, Falconer showed in \cite{Falconer} that the previous result holds for any compact subset of $\R^d$ with Hausdorff dimension larger than $\frac{d+1}{2}$. The threshold was improved to $\frac{d}{2}+ \frac{1}{3}$ by Wolff \cite{Wolff2} for $d=2$, and by Erdogan \cite{Erdogan} for $d \geq 3$. For other improved results for Falconer's distance set problem see \cite{Bourgain} and \cite{Wolff1}. The best known threshold when $d \geq 3$ is odd is $\frac{d^2}{2d-1}$. This was first achieved by Du, Guth, Ou, Wang, Wilson and Zhang \cite{DuGuthOuWangWilsonZhang} when $d=3$ and generalized to all higher dimensions by Du and Zhang \cite{DuZhang}. However, the best known threshold when restricting $d \geq 2$ to an even integer is $\frac{d}{2} + \frac{1}{4}$ obtained by Guth, Iosevich, Ou and Wang \cite{GuthIosevichOuWang} when $d=2$ and Du, Iosevich, Ou, Wang and Zhang \cite{DuIosevichOuWangZhang} for higher even dimensions.\\

There are interesting variations of the Falconer distance problem, for instance, how large does the Hausdorff dimension of a compact set $E \subset \R^d$, $d \geq 2$, need to be to ensure that the set $T_2(E)$ of all the different triangles that can be determined by triples of points of $E$ has positive three-dimensional Lebesgue measure, where different refers to triangles which are non congruent. More precisely $T_2(E) := E \times E \times E / \sim$, where $\sim$ denotes a congruence relation in $E \times E \times E$ given by
\begin{center}
  $(x^1,x^2,x^3) \sim (y^1,y^2,y^3)$ if $y^i= \tau + Rx^i$, $i=1,2,3$,
\end{center}
where $R$ is an element of the orthogonal group $O(d)$, and $\tau$ is a translation in $\R^d$. This problem was first addressed in the case $d=2$ by Greenleaf and Iosevich \cite{GreenleafIosevich}, where they showed that if $\dim_H(E)> \frac{7}{4}$, then $ \mathcal{L}^3(T_2(E))>0$. Similar results for higher dimensions as well as $k-$point configurations were obtained in \cite{GrafakosGreenleafIosevichPalsson}, \cite{GreenleafIosevichLiuPalsson1} and \cite{GreenleafIosevichLiuPalsson2}.

If a subset $E \subset \R^d$ has nonempty interior, then it has positive Lebesgue measure, however the converse is not always true. Thus, we might ask a stronger version of the Falconer distance problem, that is, for a given compact set $E \subset \R^d$, how large does its Hausdorff dimension need to be to guarantee that its distance set contains an interval. This question was initially approached by Mattila and Sj\"{o}lin in \cite{MattilaSjolin}, where they proved that if for a given compact set $E \subset \R^d$, we have that $\dim_{H}(E)>\frac{d+1}{2}$, then $\Delta(E)$ has nonempty interior. This result is nowadays known as the Mattila-Sj\"{o}lin theorem. Iosevich, Taylor and Mourgoglou   \cite{IosevichMourgoglouTaylor} not only gave an alternative proof of the Mattila-Sj\"{o}lin theorem, but they also concluded that the set $\Delta (E) = \left \{ \left \| x-y \right \|_{B} : x,y \in E \right \}$ contains an interval where $ \left \| \cdot \right \|_{B}$ denotes the metric induced by the norm defined by a bounded convex body $B$ with non-vanishing curvature. Greenleaf, Iosevich and Taylor \cite{GreenleafIosevichTaylor1} extended the Mattila-Sj\"{o}lin theorem to more general $2-$point configurations. Most recently, they upgraded their result to a general class of $k-$point configurations \cite{GreenleafIosevichTaylor2}. \\

Consider a compact set $E \subset \R^d$. The classic side-angle-side rule, used to prove whether a given set of triangles are congruent allows us to conclude that there is a one-to-one correspondence between the set $T_2(E)$ and the set
\begin{center}
$\dps \Delta_{\text{tri}}(E) = \left \{ (t,r,\alpha):  |x-z|=t, |y-z|=r \text{ and }  \alpha = \alpha(x,z,y), \  x,y,z \in E \right \} $,
\end{center}
where $\alpha (x,z,y)$ denotes the angle formed by $x$, $y$ and $z$ , centered at $z$. We will show that if $E$ has a large enough Hausdorff dimension, then the set $\Delta_{\text{tri}}(E)$ has nonempty interior. We now state our main result.

\begin{thm}
\label{thm:MainResult}
Let $E \subset \R^d$ be a compact set where $d \geq 4$. If $\dps \dim_{\textit{H}}(E) > \frac{2}{3}d+1$, then $\Delta_{\text{tri}}(E)$ has nonempty interior.
\end{thm}

In terms of sharpness we note that as the Mattila-Sj\"{o}lin type theorem for triangles implies the corresponding Falconer type theorem we inherit the sharpness examples from the Falconer setting. For triangles there is the trivial sharpness example of $\frac{d}{2}$, see for example \cite{GreenleafIosevich,GrafakosGreenleafIosevichPalsson}, which is inherited from the original Falconer distance problem, for if there are many triangles, there must be many distances. The only non-trivial sharpness example that exists for Falconer type theorems for triangles is a dimensional threshold of $\frac{3}{2}$ in the plane due to Erdo\u{g}an and Iosevich which appeared in \cite{GreenleafIosevichLiuPalsson1}. As our result is only in dimensions $4$ and higher we are left with a gap between the dimensional threshold of $\frac{2}{3}d+1$ and the trivial sharpness example up to the threshold $\frac{d}{2}$. In the setting of the original Mattila-Sj\"{o}lin theorem the situation is similar, the only sharpness examples that are known are those that come from the Falconer distance problem.

With the lack of sharpness examples for the Mattila-Sj\"{o}lin type theorems that differ from those that arise for the Falconer type theorems, one might wonder if in fact one should expect the same dimensional thresholds for the two questions. Despite lack of evidence in the Euclidean setting, there is evidence in the finite field setting that the thresholds may differ. Murphy, Petridis, Pham, Rudnev and Stevens \cite{MurphyPetridisPhamRudnevStevens} proved that if $E \subset \mathbb{F}^2_q$, $q$ prime, and the cardinality of $E$ is greater than $Cq^{\frac{5}{4}}$, that is $|E| \geq Cq^{\frac{5}{4}}$, for some $C>0$, then its distance set $\Delta(E)$ contains a positive proportion of $\mathbb{F}_q$. This result can be thought of as an analogue to a Falconer type result. Murphy and Petridis \cite{MurphyPetridis} proved that if $E \subset \mathbb{F}^2_q$ and $|E| \approx q^{\frac{4}{3}}$, then one cannot, in general, expect the distance set to contains all $\mathbb{F}_q$, which can be thought of as the impossibility of a Mattila-Sj\"{o}lin type result. Yet Iosevich and Rudnev \cite{IosevichRudnev} proved that if $E \subset \mathbb{F}^d_q$, $d \geq 2$, such that $|E| \gtrsim q^{\frac{d+1}{2}}$, then its distance set $ \Delta(E)$ contains all of $\mathbb{F}_q$ which establishes a Mattila-Sj\"{o}lin type result at this higher threshold. 

We now return to some consequences of our main theorem. Let $\mu$ be a Frostman probability measure supported on $E$, so $\mu(B(x,r)) \leq c_{\mu}r^s$ for some $s>0$ and for all $x \in \R^d$, $r>0$. By Frostman's Lemma the parameter $s$ can be taken all the way up to (but not including) $\dim_{\textit{H}}(E)$. We can now define a measure $\delta(\mu)$ on $ \Delta_{\text{tri}}(E)$ by setting, for any continuous function $h$,
\be \label{eqn:defi}
  \int h(t,r,\alpha) \df \delta(\mu)(t,r,\alpha) = \iiint h(|x-z|,|y-z|, \alpha(x,z,y))  \df \mu (x) \df \mu(y) \df \mu (z).
\ee
Thus $\delta(\mu)$ is the image of $\mu \times \mu \times \mu$ under the map $(x,y,z) \to (|x-z|,|y-z|, \alpha(x,z,y))$. For a smooth compactly supported function $f$, $\delta(f)$ is also a function given by
\be \label{eqn:eq1}
    \delta(f)(t,r,\alpha) = \iint f(z+tx)f(z+ry)f(z)t^{d-1}r^{d-1} \df \sigma_{\alpha}(x,y) \df z
\ee
which resembles the associated spherical means used in \cite{Liu}. The measure $\df \sigma_{\alpha}$ given in the right hand side of (\ref{eqn:eq1}) denotes the normalized surface area measure over the surface
\begin{equation*}
    \lla{ (x,y) \in \R^{2d}: |x|=|y|=1 \text{\ and \ } y=g_{\alpha}x },
\end{equation*}
where $g_{\alpha} \in O(d)$ denotes some rotation by the angle $\alpha$. To show this note
\begin{multline*}
    \int g(t,r,\alpha) \df \delta (f)(t,r,\alpha) = \iiint g(|x-z|,|y-z|, \alpha(x,z,y)) \\ f(x) f(y) f(z) \df \mu (x) \df \mu(y) \df \mu (z),
\end{multline*}
for any continuous function $g$ with compact support. By a change of variables we have
\begin{multline*}
\int g(t,r,\alpha) \df \delta (f)(t,r,\alpha) = \iiint g(|x|,|y|,\alpha(z+x,z,z+y)) \\
f(z+x) f(z+y) f(z) \df x \df y \df z
\end{multline*}
and thus, by using polar coordinates we obtain
\begin{multline*}
\int g(t,r,\alpha) \df \delta (f)(t,r,\alpha)= \iiint g(t,r,\alpha) \\ \pa{ \iint f(z+tx) f(z+ry) f(z) t^{d-1} r^{d-1} \df \sigma_{\alpha}(x,y) \df z }  \df t \df r \df \alpha.
\end{multline*}
The authors in \cite{HarangiKeletiKissMagaMatheMattilaStrenner} showed that if a compact set $E \subset \mathbb{R}^d$, $d \geq 2$, has a Hausdorff dimension greater than $d-1$, then any angle $\alpha \in [0,\pi]$ can be formed by triples of points of $E$. Therefore, the set of angles formed by triples of points of $E$ has nonempty interior. The following corollary, which is an immediate consequence of our main result, establishes a related result.
\begin{cor}
  Let $E \subset \mathbb{R}^d$ be a compact set where $d \geq 4$. If $\displaystyle \dim_{\textit{H}}(E) > \frac{2}{3}d+1$, then the angles set
  \begin{center}
    $\displaystyle \mathcal{A}(E) = \left \{ \alpha(x,z,y): x,y,z \in E \text{\ distinct } \right \},$
  \end{center}
   where $\alpha (x,z,y)$ denotes the angle formed by $x$, $y$ and $z$, centered at $z$, has nonempty interior.
\end{cor}
Thus for sets of dimension greater than $\frac{2}{3} d + 1$ one can guarantee lots of angles, we emphasize though that we can't say which angles are guaranteed to arise which is weaker than the statements in \cite{HarangiKeletiKissMagaMatheMattilaStrenner}. For Falconer type theorems on the Lebesgue measure of the angle set $\dps \mathcal{A}(E)$ being positive see \cite{IosevichMourgoglouPalsson,IosevichPalsson}.\\

For completeness we point out the following incidence estimate as a consequence of our main result.

\begin{cor}\label{incthm}
    Let $E \subset \mathbb{R}^d$ be a compact set where $d \geq 4$ and $\mu$ a Frostman probability measure supported on $E$, meaning $\mu(B(x,R)) \leq c_{\mu}R^s$ for all $x \in \mathbb{R}^d$, $R>0$. If $\displaystyle \dim_{\textit{H}}(E) > \frac{2}{3}d+1$, then for all $0<\varepsilon < t,r,\alpha$ such that $\alpha +\varepsilon < \frac{\pi}{2}$,
    \begin{multline*}
    \mu \times \mu \times \mu    \left ( \left \{ (x,y,z): t \leq |x-z| \leq t + \varepsilon ,\ r \leq |y-z| \leq r + \varepsilon, \right. \right. \\
    \left. \left.  \alpha \leq \alpha(x,y,z) \leq \alpha + \varepsilon  \right \} \right )  \leq \widetilde{C}(t,r,\alpha) \varepsilon^3.
    \end{multline*}
\end{cor}
Stronger results were obtained in \cite{GrafakosGreenleafIosevichPalsson,GreenleafIosevich}. For sharpness in incidence estimates see \cite{DeWittFordGoldsteinMillerMorelandPalssonSenger} and the references contained therein.

\subsection{Contributions of this work:} As a consequence of the result given by the authors in \cite{GreenleafIosevichTaylor2}, we have that for a compact set $E \subset \R^d$, if $\dim_{H}(E)> \frac{2d+3+\beta}{3}$, then the $3-$point $\Phi$-configuration,
\benn
  \tri_{\Phi}(E):= \lla{\phi(x_1,x_2,x_3) : \ x_1,x_2,x_3 \in E} \subset \R^3,
\eenn
has nonempty interior. Here $\beta \geq 0$, and $\Phi: E \times E \times E \to T \subset \R^3$ must satisfy some given conditions. It is not easy to see if there is a map $\Phi : E \times E \times E \to \Delta_{\text{tri}}(E)$ that satisfies all the conditions given in \cite{GreenleafIosevichTaylor2}, but even if we were able to define such a map the threshold obtained would not be smaller than $\frac{2}{3}d+1$. In fact, $\frac{2}{3}d+1$ is the best threshold we could expect from \cite{GreenleafIosevichTaylor2}.

\subsection{Overview of result: } We build on the techniques used by Iosevich and Liu in \cite{IosevichLiu}, where they proved the existence of a value $s_{d,\mu}<d$, such that for a compact set $E \subset \R^d$, $d \geq 4$, with $\dim_{H}(E) > s_{d,\mu} $, the set $E$ contains $3-$point configurations that form an equilateral triangle. Recently, Iosevich and Maygar \cite{IosevichMagyar} showed a similar result for simplices.\\

For the proof of our main result consider $\ve>0$ and define $\dps \mu_{\ve} := \mu \ast \phi_{\ve}$, where $\dps \phi_{\ve}= \ve^{-d} \phi \pa{\frac{.}{\ve}} $ and $\dps \phi \in C^{\infty}_0$ is supported in the unit ball with $\int \phi =1$. The proof will be done in the following steps:
\begin{itemize}
    \item {\bf Step 1:} Due to we have that $\delta(\mu_{\ve}) \to \delta(\mu)$ weakly as long as $\mu_{\ve} \to \mu$ weakly, we only need to show that $\delta(\mu_{\ve})$ converges (strongly) to a continuous function. This will be done  by following the general approach of \cite{IosevichMourgoglouTaylor}, that is, we will express $\delta(\mu_{\ve})(t,r,\alpha)$ as
\benn
   \delta(\mu_{\ve})(t,r,\alpha)= M(t,r,\alpha)+R_{\ve}(t,r,\alpha),
\eenn
    \item {\bf Step 2:} We will use Proposition \ref{Claim} given below, to show that $\dps M(t,r,\alpha)$ is a continuous function and $\dps \lim_{\ve \to 0} R_{\ve}(t,r,\alpha)=0$
\end{itemize}
This will allow us to conclude the proof of our theorem.

\subsection{Acknowledgements: }
We would like to express our deep gratitude to the anonymous referee whose suggestions have greatly contributed to this paper. 


\section{Proof of Theorem \ref{thm:MainResult}}

For our purpose it is convenient to ensure that we can extract three suitable subsets of $E$ disjoint from each other. To accomplish this stated goal, we have the following result which will be proved in section $6$.

\begin{lem}\label{PHprinc}
  Let $\mu$ be a Frostman probability measure on $E \subset \R^d$, with Hausdorff dimension greater than $\dps \frac{2}{3} d+1$, then there are positive constants $c_1$, $c_2$, $c_3$, $c_4$ and $E_1$, $E_2$, $E_3$ subsets of $E$, such that
  \begin{itemize}
    \item [(i)] $\dps \mu(E_i) \geq c_1 >0$, for $i=1,2,3$.
    \item [(ii)] $\dps \underset{1 \leq k \leq d}{\max} \left \{  \inf \left \{ | x_k - y_k |: x \in E_i, y \in E_j \text{ \ and \ } i \neq j \right \}  \right \} \geq c_2 >0 $, for $i,j=1,2,3 $.
    \item [(iii)] $\dps 0 < c_3 \leq \alpha(x,z,y) \leq c_4 < \frac{\pi}{2} $, for $x \in E_1$, $y \in E_2$ and $z \in E_3$.
  \end{itemize}
\end{lem}

Let $E_1$, $E_2$ and $E_3$ be subsets of $E$ as in Lemma \ref{PHprinc}, and let $\mu_1$, $\mu_2$ and $\mu_3$ be restrictions of $\mu$ to the sets $E_1$, $E_2$ and $E_3$, respectively. The standard proof is given in terms of $E$ and $\mu$, but it is clear that the proof also works for $E_i$ and $\mu_i$, $i=1,2,3$.\\
\begin{rek}
    We are interested in non-degenerate triangles, so Lemma \ref{PHprinc} allow us to ensure that the set of $t's$, $r's$ and $\alpha 's$ can be bounded above and below by positive constants. The reason why we need the set of $\alpha 's$ to be bounded away from $\frac{\pi}{2}$ is because we do not know how to estimate the decay of the Fourier transform of the surface measure $\sigma_{\alpha}$ (See proof of Lemma \ref{lem2.2}).
\end{rek}
\subsection{Proof of Step 1:} From (\ref{eqn:eq1}) we have
\benn
\delta(\mu_{\ve})(t,r,\alpha) = \iint \mu_{\ve}(z+tx)  \mu_{\ve}(z+ry) \mu_{\ve}(z) t^{d-1} r^{d-1} \df \sigma_{\alpha}(x,y) \df z.
\eenn
By the Fourier inversion formula
\begin{eqnarray*}
\delta(\mu_{\ve})(t,r,\alpha) & = & \iiiint \h{\mu}_{\ve}(\zeta) e^{2 \pi i \zeta \cdot z} \h{\mu}_{\ve}(\xi) e^{2 \pi i \xi \cdot (z+tx)} \h{\mu}_{\ve}(\eta) e^{2 \pi i \eta \cdot (z+ry)} t^{d-1} r^{d-1}  \df \sigma_{\alpha}(x,y) \df \xi \df \zeta d\eta \\
                         & = & \iiiint \h{\mu}_{\ve}(\zeta) \h{\mu}_{\ve}(\xi) \h{\mu}_{\ve}(\eta) e^{2 \pi i (z,z,z) \cdot (\xi ,\eta ,\zeta)} e^{2 \pi i (x,y) \cdot (t \xi , r \eta)} t^{d-1} r^{d-1} \df \sigma_{\alpha}(x,y) \df \xi \df \zeta \df \eta \\
                         & = & \iiint \h{\mu}_{\ve}(\zeta) \h{\mu}_{\ve}(\xi) \h{\mu}_{\ve}(\eta) e^{2 \pi i (z,z,z) \cdot (\xi ,\eta ,\zeta)} \h{\sigma}_{\alpha}(-t \xi, -r \eta) t^{d-1} r^{d-1} \df \xi \df \zeta \df \eta \\
\end{eqnarray*}
due to
\begin{center}
    $\dps e^{2 \pi i (z,z,z) \cdot (\xi ,\eta ,\zeta)} = \delta(\xi+\eta+\zeta)$
\end{center}
we have
\begin{eqnarray*}
\delta(\mu_{\ve})(t,r,\alpha) & = & \iint \h{\mu}_{\ve}(-\xi -\eta ) \h{\mu}_{\ve}(\xi) \h{\mu}_{\ve}(\eta) \h{\sigma}_{\alpha}(- t\xi, - r\eta) t^{d-1} r^{d-1} \df \xi \df \eta, \\
\end{eqnarray*}
and because of $\dps \mu_{\delta} = \mu \ast \phi_{\delta}$, the last equation can be expressed as
\begin{multline*}
\delta(\mu_{\ve})(t,r,\alpha)  = \iint \h{\mu}(-\xi -\eta) \h{\mu}(\xi) \h{\mu}(\eta) \h{\sigma}_{\alpha}(- t\xi, - r\eta) t^{d-1} r^{d-1} \df \xi \df \eta \\
- \iint \h{\mu}(-\xi -\eta) \corch{1-\h{\phi}_{\ve} (-\xi -\eta) \h{\phi}_{\ve}(\xi) \h{\phi}_{\ve}(\eta) } \h{\mu}(\xi) \h{\mu}(\eta) \h{\sigma}_{\alpha}(- t\xi,-r\eta) t^{d-1} r^{d-1}  \df \xi \df \eta.
\end{multline*}
Let
\begin{equation*}
    M(t,r,\alpha) := \iint \h{\mu}(-\xi -\eta) \h{\mu}(\xi) \h{\mu}(\eta) \h{\sigma}_{\alpha}(-t\xi, -r\eta) t^{d-1} r^{d-1} \df \xi \df \eta
\end{equation*}
and
\begin{multline*}
    R_{\ve}(t,r,\alpha) := - \iint \h{\mu}(-\xi -\eta) \left [ 1 - \h{\phi}_{\ve} (-\xi-\eta) \h{\phi}_{\ve}(\xi) \h{\phi}_{\ve}(\eta) \right ] \\ \h{\mu}(\xi) \h{\mu}(\eta)  \h{\sigma}_{\alpha}(-t\xi,-r\eta) t^{d-1} r^{d-1} \df \xi \df \eta.
\end{multline*}
\subsection{Proof of Step 2:} We will need the following
\begin{prop}\label{Claim}
For $R>0$ we have
   \begin{equation*}
  \iint \limits_{ |\xi|, |\eta| > R } |\h{\mu}(\xi +\eta)| |\h{\mu}(\xi)||\h{\mu}(\eta)||\h{\sigma}_{\alpha}(t\xi,r\eta)| t^{d-1} r^{d-1}  \df \xi \df \eta \lesssim \widehat{C}(t,r,\alpha) R^{-\frac{3s-2d-3}{2}}
  \end{equation*}
where $\dps \widehat{C}(t,r,\alpha) = C_1(\alpha) (max \lla{t,r })^{-\frac{1}{2}} t^{\frac{d}{2}} r^{\frac{d}{2}} $, and $C_1(\alpha)$ is an integrable function.
\end{prop}
Proof of Proposition \ref{Claim} can be found in Section $2$. To show that $R_{\ve}(t,r,\alpha) \to 0$ as $\ve \to 0$, note
\begin{eqnarray*}
  |R_{\ve}(t,r,\alpha)| & \leq & \iint \limits_{ |\xi|, |\eta| > \frac{1}{2\ve}}|\h{\mu}(\xi +\eta)| |\h{\mu}(\xi)|| \h{\mu}(\eta)||\h{\sigma}_{\alpha}( t\xi, r\eta)| t^{d-1} r^{d-1}  \df \xi \df \eta .
\end{eqnarray*}
Therefore, by Proposition \ref{Claim} we have
\begin{equation*}
  |R_{\ve}(t,r,\alpha)| \lesssim \widehat{C}(t,r,\alpha) (2\ve)^{\frac{3s-2d-3}{2}}.
\end{equation*}
Then $\dps \lim_{\ve \to 0} R_{\ve}(t,r,\alpha)=0$, if $\dps s>\frac{2d}{3}+1$. \\

To show that $M(t,r,\alpha)$ is continuous consider $R>0$ and note that if $|\xi|,|\eta|>R$, then by Lemma \ref{lem2.2} (See proof of Proposition \ref{Claim} below) the integrand of $M(t,r,\alpha)$ is bounded by an $L^1$ function. If $|\xi|,|\eta| \leq R$, then the compactness of $S^{2d-1}$ and the fact that $\dps \h{\sigma}_{\alpha}(- t\xi,-r\eta)$ is given in terms of local coordinates, we have that the integrand of $M(t,r,\alpha)$ is also be bounded by an $L^1$ function.\\
 Moreover, due to $\dps\h{\sigma}_{\alpha}(- t\xi, - r\eta) = \int e^{-2 \pi i U_{\alpha} \cdot (- t\xi, - r\eta) } \df u_1 \df u' \df v'$ , where $U_{\alpha}$ is a local parametrization (see proof of Lemma 2.1 below) given by
 \begin{align*}
      U_{\alpha} = & \left ( u_1, u', 1-\frac{|u|^2}{2} ; \sin(\alpha)+ \cos(\alpha)u_1 + \cot(\alpha) (u' \cdot v') - \frac{1}{2\sin(\alpha)} |v'|^2 - \frac{\cos^2(\alpha)}{2\sin(\alpha)} |u'|^2  \right. \\
     & \left. - \frac{\sin(\alpha)}{2} u^2_1, v', \cos(\alpha) - \sin(\alpha)u_1 - u' \cdot v' - \frac{\cos(\alpha)}{2} |u'|^2  - \frac{\cos(\alpha)}{2} u^2_1  \right ),
\end{align*}
 the integrand of $M(t+h_1,r+h_2,\alpha+h_3)$ converges to the integrand of $M(t,r,\alpha)$ as $h_1$, $h_2$ and $h_3$ tend to zero. Therefore, by the Dominated Convergence Theorem it follows that $M(t,r,\alpha)$ is continuous.

\section{Proof of Proposition \ref{Claim}}

For the first part of the proof we follow the ideas given in \cite{IosevichLiu}. Let $R$ be any positive number. We would like to show
  \begin{equation}\label{eqn:eq2}
  \iint \limits_{ |\xi|, |\eta|> R } |\h{\mu}(\xi +\eta)| |\h{\mu}(\xi)|| \h{\mu}(\eta)||\h{\sigma}_{\alpha}( t\xi, r\eta)| t^{d-1} r^{d-1}  \df \xi \df \eta \lesssim \widehat{C}(t,r,\alpha) R^{-\frac{(3s-2d-3)}{2}}.
  \end{equation}
Where $\widehat{C}(t,r,\alpha) = C_1(\alpha) (max \lla{ t,r})^{-\frac{1}{2}} t^{\frac{d}{2}} r^{\frac{d}{2}} $. To show (\ref{eqn:eq2}) we will need the following
\begin{lem}\label{lem2.2}
Consider $0 < \alpha < \frac{\pi}{2}$, then
  \benn
    |\h{\sigma}_{\alpha}(\xi,\eta)| \leq C_1(\alpha) |\xi + g_{\alpha}\eta|^{-\frac{1}{2}} |\xi|^{-\frac{(d-2)}{2}} |\eta|^{-\frac{(d-2)}{2}} (\sin( < \xi , \eta >)) ^{-\frac{(d-2)}{2}}
  \eenn
  Where $g_{\alpha} \in O(d)$ is some rotation by $\alpha$, $C_1(\alpha)$ is an integrable function, and $<\xi,\eta>$ denotes the angle between $\xi$ and $\eta$.
\end{lem}
In \cite{IosevichLiu} the authors proved Lemma \ref{lem2.2} for the case $\alpha=\frac{\pi}{3}$, their proof can be easily generalized to any $0 < \alpha < \frac{\pi}{2} $. For the sake of completeness the proof of Lemma \ref{lem2.2} can be found in Section $4$. Note that if $|\xi|+|\eta| \approx R2^{j}$, then at least two of $|\xi|$, $|\eta|$, $|\xi+\eta|$ are $\approx R2^j$. Thus we have the following cases,\\
{\bf Case 1:} $\dps |\xi| \approx |\eta| \approx R2^j$.\\
By lemma \ref{lem2.2}
\begin{multline*}
 \iint \limits_{|\xi| \approx |\eta| \approx R2^j} |\h{\mu}(\xi)||\h{\mu}(\eta)||\h{\mu}(\xi + \eta)||\h{\sigma}_{\alpha}(t\xi,r\eta)| \df \xi \df \eta \\
\lesssim  C(t,r,\alpha) R^{2-d} 2^{-j(d-2)} \iint \limits_{|\xi| \approx |\eta| \approx R2^j} |\h{\mu}(\xi)||\h{\mu}(\eta)||\h{\mu}(\xi + \eta)||\xi + g_{\alpha}\eta|^{-\frac{1}{2}} (\sin( <\xi,\eta > )) ^{-\frac{(d-2)}{2}} \df \xi \df \eta,
\end{multline*}
where $\dps C(t,r,\alpha) = C_1(\alpha)(max \lla{t,r})^{-\frac{1}{2}} t^{-\frac{d-2}{2}} r^{-\frac{d-2}{2}} $. Fix $\eta$
\begin{multline*}
\int \limits_{|\xi| \approx R2^j} |\h{\mu}(\xi+\eta)| |\xi + g_{\alpha}\eta|^{-\frac{1}{2}} (\sin( <\xi,\eta>)) ^{-\frac{(d-2)}{2}} \df \xi \\ \lesssim \pa{\int \limits_{ |\xi| \approx R2^j} |\h{\mu}(\xi+\eta)|^2 \df \xi}^{\frac{1}{2}} \pa{ \int \limits_{|\xi| \approx R2^j} |\xi + g_{\alpha}\eta|^{-1} (\sin( < \xi , \eta > )) ^{-(d-2)} \df \xi }^{\frac{1}{2}}
\end{multline*}
since $\mu(B(x,r)) \lesssim r^s $, then $\dps \int \limits_{ |\xi| \approx R} |\h{\mu}(\xi)|^2 \lesssim R^{d-s}$. Thus we have
\begin{multline*}
\int \limits_{|\xi| \approx R2^j} |\h{\mu}(\xi+\eta)| |\xi + g_{\alpha}\eta|^{-\frac{1}{2}} (\sin( < \xi , \eta > )) ^{-\frac{(d-2)}{2}} \df \xi \\ \lesssim \pa{R2^j }^{\frac{d-s}{2}} \pa{ \int \limits_{|\xi| \approx R2^j} |\xi + g_{\alpha}\eta|^{-1} (\sin( < \xi , \eta > )) ^{-(d-2)} \df \xi }^{\frac{1}{2}}
\end{multline*}
The following lemma is also proved in \cite{IosevichLiu} for the case $\displaystyle \alpha = \frac{\pi}{3}$ and it can be generalized for our purpose. The proof can be found in Section $5$.
\begin{lem} \
\benn \label{lem2.3}
     \int \limits_{|\xi| \approx R2^j} |\xi + g_{\alpha}\eta|^{-1} (\sin( < \xi , \eta > ))^{-(d-2)} \df \xi \lesssim  \pa{ R2^j }^{-1+d}.
\eenn
\end{lem}
Thus by Lemma \ref{lem2.3}
\benn
  \int \limits_{|\xi| \approx R2^j} |\h{\mu}(\xi+\eta)| |\xi + g_{\alpha}\eta|^{-\frac{1}{2}} (\sin( < \xi , \eta >)) ^{-\frac{(d-2)}{2}} \df \xi \lesssim \pa{ R2^j }^{\frac{2d-s-1}{2}}.
\eenn
Similarly, for $\xi$ fixed we have
\benn
\int \limits_{|\eta| \approx R2^j} |\h{\mu}(\xi+\eta)| |\xi + g_{\alpha}\eta|^{-\frac{1}{2}} (\sin( < \xi , \eta >)) ^{-\frac{(d-2)}{2}} \df \eta \lesssim \pa{R2^j }^{\frac{2d-s-1}{2}}.
\eenn
By Shur's Test (see \cite{Grafakos})
\begin{eqnarray*}
& & \iint \limits_{|\xi| \approx |\eta| \approx R2^j} |\h{\mu}(\xi)||\h{\mu}(\eta)||\h{\mu}(\xi + \eta)||\h{\sigma}_{\alpha}(t\xi,r\eta)| \df \xi \df \eta \\
& \lesssim & C(t,r,\alpha) R^{2-d} 2^{-j(d-2)} \pa{R2^j}^{\frac{2d-s-1}{2}} \pa{ \int \limits_{|\xi| \approx R2^j} |\h{\mu}(\xi)|^2 \df \xi } \\
& \lesssim &  C(t,r,\alpha) R^{- \pa{\frac{3s-2d-3}{2}} } 2^{-j  \pa{\frac{3s-2d-3}{2}} } \\
\end{eqnarray*}
and thus
\begin{multline*}
  \iint \limits_{|\xi| \approx |\eta| \approx R2^j} |\h{\mu}(\xi)||\h{\mu}(\eta)||\h{\mu}(\xi + \eta)||\h{\sigma}_{\alpha}(t\xi,r\eta)| t^{d-1}r^{d-1} \df \xi \df \eta \\ \lesssim t^{d-1}r^{d-1} C(t,r,\alpha) R^{-\frac{3s-2d-3}{2}} 2^{-j  \pa{ \frac{3s-2d-3}{2}} }.
\end{multline*}
{\bf Case 2:} $\dps|\xi| \approx |\xi + \eta| \approx R2^j$.\\
In \cite{IosevichLiu} the authors used the fact that $\dps \alpha =\frac{\pi}{3}$ to reduce everything to case 1. We will decompose the domain of integration to find an upper bound. By Lemma \ref{lem2.2} we have
\begin{multline*}
 \iint \limits_{|\xi| \approx |\xi +\eta| \approx R2^j} |\h{\mu}(\xi)||\h{\mu}(\eta)||\h{\mu}(\xi + \eta)||\h{\sigma}_{\alpha}(t\xi,r\eta)| \df \xi \df \eta  \\
\leq C_1(\alpha) \iint \limits_{|\xi| \approx |\xi+\eta| \approx R2^j} |\h{\mu}(\xi)||\h{\mu}(\eta)||\h{\mu}(\xi + \eta)||t\xi + g_{\alpha}r\eta|^{-\frac{1}{2}} |t\xi|^{-\frac{(d-2)}{2}} |r\eta|^{-\frac{(d-2)}{2}} (\sin( < \xi , \eta > ))^{-\frac{(d-2)}{2}} \df \xi \df \eta.
\end{multline*}
The last expression is bounded by
\begin{multline} \label{eqn:eq3}
  C(t,r,\alpha) R^{\frac{(2-d)}{2}} 2^{-j\frac{(2-d)}{2}} \iint \limits_{|\xi| \approx |\xi+ \eta| \approx R2^j} |\h{\mu}(\xi)||\h{\mu}(\eta)||\h{\mu}(\xi + \eta)||\xi + g_{\alpha}\eta|^{-\frac{1}{2}} \\ |\eta|^{-\frac{(d-2)}{2}} (\sin( < \xi , \eta >)) ^{-\frac{(d-2)}{2}} \df \xi \df \eta
\end{multline}
Where $\displaystyle C(t,r,\alpha) = C_1(\alpha) (max\left \{ t,r \right \})^{-\frac{1}{2}} t^{-\frac{d-2}{2}} r^{-\frac{d-2}{2}} $. Note that the integral above is equal to
\begin{multline*}
  \iint \limits_{ \begin{matrix}  |\xi| \approx |\xi+ \eta| \approx R2^j \\ \frac{R2^j}{2} \leq |\eta| \leq 2R2^j \end{matrix} } |\h{\mu}(\xi)||\h{\mu}(\eta)||\h{\mu}(\xi + \eta)||\xi + g_{\alpha}\eta|^{-\frac{1}{2}} |\eta|^{-\frac{(d-2)}{2}}(\sin( < \xi , \eta >)) ^{-\frac{(d-2)}{2}} \df \xi \df \eta  \\
  + \sum^{j+1}_{k=3} \iint \limits_{ \begin{matrix} |\xi| \approx |\xi+ \eta| \approx R2^j \\ \frac{2R2^j}{2^k} \leq |\eta| \leq \frac{2R2^j}{2^{k-1}} \end{matrix} } |\h{\mu}(\xi)||\h{\mu}(\eta)||\h{\mu}(\xi + \eta)||\xi + g_{\alpha}\eta|^{-\frac{1}{2}} |\eta|^{-\frac{(d-2)}{2}}(\sin( < \xi , \eta >))^{-\frac{(d-2)}{2}} \df \xi \df\eta.
\end{multline*}
  The first term of the sum above can be reduce to the Case 1: $\dps |\xi| \approx |\eta| \approx R2^j$, that is,
\begin{multline*}
    \iint \limits_{ \begin{matrix} |\xi| \approx |\xi+ \eta| \approx R2^j \\ \frac{R2^j}{2} \leq |\eta| \leq 2R2^j \end{matrix} } |\h{\mu}(\xi)||\h{\mu}(\eta)||\h{\mu}(\xi + \eta)||\xi + g_{\alpha}\eta|^{-\frac{1}{2}} |\eta|^{-\frac{(d-2)}{2}}(\sin( < \xi , \eta >)) ^{-\frac{(d-2)}{2}} \df \xi \df \eta \\ \lesssim C(t,r,\alpha) R^{(2-d)} 2^{-j(d-2)} (R2^j)^{\frac{2d-s-1}{2}}(R2^j)^{d-s}.
\end{multline*}
  For the other terms of the sum note that for each $k$, and $\xi$ fixed we have
\begin{align*}
  & \int \limits_{ \begin{matrix} |\xi+ \eta| \approx R2^j \\ \frac{2R2^j}{2^k} \leq |\eta| \leq \frac{2R2^j}{2^{k-1}} \end{matrix} }  |\h{\mu}(\xi+\eta)| |\xi + g_{\alpha}\eta|^{-\frac{1}{2}} |\eta|^{-\frac{d-2}{2}} (\sin( < \xi , \eta > )) ^{-\frac{(d-2)}{2}} \df \eta \\
  & \lesssim  \pa{ \int \limits_{  |\eta| \approx \frac{R2^j}{2^{k-1}} }  |\h{\mu}(\xi+\eta)|^2 \df \eta }^{\frac{1}{2}} \pa{ \int \limits_{ \begin{matrix} |\xi+ \eta| \approx R2^j \\ \frac{2R2^j}{2^k} \leq |\eta| \leq \frac{2R2^j}{2^{k-1}} \end{matrix} } |\xi + g_{\alpha}\eta|^{-1}  |\eta|^{-(d-2)} (\sin( < \xi , \eta > )) ^{-(d-2)} \df \eta}^{\frac{1}{2}} \\
  & \lesssim \pa{\frac{R2^j}{2^{k-1}} }^{\frac{2-s}{2}} \pa{ R2^j }^{-\frac{1}{2}} \pa{\frac{R2^j}{2^{k-1}}}^{\frac{d}{2}}
 \end{align*}
where the last inequality follows from Lemma $\ref{lem2.3}$. Similarly for $\eta$ fixed
\begin{multline*}
   \int \limits_{ \begin{matrix} |\xi+ \eta| \approx R2^j \\ \frac{2R2^j}{2^k} \leq |\eta| \leq \frac{2R2^j}{2^{k-1}} \end{matrix} }  |\h{\mu}(\xi+\eta)| |\xi + g_{\alpha}\eta|^{-\frac{1}{2}} |\eta|^{- \frac{d-2}{2}} (\sin( < \xi , \eta > )) ^{-\frac{(d-2)}{2}} \df \xi \\ \lesssim \pa{R2^j}^{\frac{2d-s-1}{2}} \pa{\frac{R2^j}{2^{k-1}}}^{-\frac{(d-2)}{2}}.
 \end{multline*}
Thus by Shur's test we have
\begin{multline*}
  \iint \limits_{ \begin{matrix} |\xi| \approx |\xi+ \eta| \approx R2^j \\ \frac{2R2^j}{2^k} \leq |\eta| \leq \frac{2R2^j}{2^{k-1}} \end{matrix} } |\h{\mu}(\xi)||\h{\mu}(\eta)||\h{\mu}(\xi)||\h{\mu}(\eta)| |\h{\mu}(\xi+\eta)| |\xi + g_{\alpha}\eta|^{-\frac{1}{2}} |\eta|^{- \frac{d-2}{2}} \\ (\sin(< \xi , \eta >)) ^{-\frac{(d-2)}{2}} \df \xi \df \eta \lesssim   \pa{R2^j}^{\frac{d-s+1}{2}} \pa{\frac{1}{2^{k-1}}}^{\frac{4-s}{4}} \pa{\frac{R2^j}{2^{k-1}} }^{d-s}.
 \end{multline*}
Therefore (\ref{eqn:eq3}) is bounded by
\be \label{eqn:c}
  C(t,r,\alpha) R^{\frac{3+2d-3s}{2}} \corch{ 2^{-j \pa{\frac{3s-2d-3}{2}} }   +  2^{-j \pa{\frac{s-2}{4}} } \sum^{j-2}_{k=0} 2^{k \pa{ \frac{4+4d-5s}{4} }}}.
\ee
Note that $\dps \frac{2d}{3}+1 < \frac{4+4d}{5}$. Thus, if $\dps\frac{2d}{3}+1 < s < \frac{4+4d}{5}$ we have
  \begin{eqnarray*}
 (\ref{eqn:c}) & \leq & C(t,r,\alpha) R^{\frac{3+2d-3s}{2}} \left [  2^{-j \left ( \frac{3s-2d-3}{2} \right ) }   +  2^{-j \left ( \frac{s-2}{4} \right ) } \sum^{j-2}_{k=0} 2^{j \left ( \frac{4+4d-5s}{4} \right ) }  \right ] \\
 & = & C(t,r,\alpha) R^{\frac{3+2d-3s}{2}} j 2^{-j \left ( \frac{3s-2d-3}{2} \right ) } \\
\end{eqnarray*}
and thus
\begin{multline*}
  \iint \limits_{|\xi| \approx |\xi+ \eta| \approx R2^j} |\h{\mu}(\xi)||\h{\mu}(\eta)||\h{\mu}(\xi + \eta)||\h{\sigma}_{\alpha}(t\xi,r\eta)| t^{d-1}r^{d-1} \df \xi \df \eta \\ \lesssim t^{d-1}r^{d-1} C(t,r,\alpha) R^{-\frac{3s-2d-3}{2}} j 2^{-j \left ( \frac{3s-2d-3}{2} \right ) }
\end{multline*}
If $\dps \frac{4+4d}{5} <  s $ we have
  \begin{eqnarray*}
 (\ref{eqn:c}) & \leq & C(t,r,\alpha) R^{\frac{3+2d-3s}{2}} \left [  2^{-j \left ( \frac{3s-2d-3}{2} \right ) }   +  (j-1)2^{-j \left ( \frac{s-2}{4} \right ) }  \right ] \\
 & = & C(t,r,\alpha) R^{\frac{3+2d-3s}{2}} j 2^{-j \left ( \frac{s-2}{4} \right ) } \\
\end{eqnarray*}
and thus
\begin{multline*}
  \iint \limits_{|\xi| \approx |\xi+ \eta| \approx R2^j} |\h{\mu}(\xi)||\h{\mu}(\eta)||\h{\mu}(\xi + \eta)||\h{\sigma}_{\alpha}(t\xi,r\eta)| t^{d-1}r^{d-1} \df \xi \df \eta \\ \lesssim t^{d-1}r^{d-1} C(t,r,\alpha) R^{-\frac{3s-2d-3}{2}} j 2^{-j \pa{\frac{s-2}{4}} }
\end{multline*}
{\bf Case 3:} $\dps |\eta| \approx |\xi + \eta| \approx R2^j$. This case is similar to case 2.\\

Therefore, in any of the previous cases we have
\begin{eqnarray*}
   & & \iint \limits_{ |\xi|, |\eta|> R} |\h{\mu}(\xi)||\h{\mu}(\eta)||\h{\mu}(\xi + \eta)||\h{\sigma}_{\alpha}(t\xi,r\eta)| t^{d-1} r^{d-1} \df \xi \df \eta \\
  & = & t^{d-1} r^{d-1} \sum_{j} \iint \limits_{ |\xi|+|\eta| \approx R2^j} |\h{\mu}(\xi)||\h{\mu}(\eta)||\h{\mu}(\xi + \eta)||\h{\sigma}_{\alpha}(t\xi,r\eta)|  \df \xi \df \eta \\
  & \lesssim & t^{d-1} r^{d-1} C(t,r,\alpha) R^{-\frac{3s-2d-3}{2}},
\end{eqnarray*}
if $\dps s>\frac{2d}{3}+1$.


\section{Proof of Lemma \ref{lem2.2}}

As we mentioned above, the proof of this lemma can be found in \cite{IosevichLiu}, where the authors worked the case $\alpha=\frac{\pi}{3}$. They followed the ideas given in \cite{GrafakosGreenleafIosevichPalsson}, where it was proved that $|\h{\sigma}(\xi,\eta)| \lesssim (1+|\xi|+|\eta|)^{-\frac{(d-1)}{2}}$. Throughout the proof we will denote $x = (x_1,x',x_d) \in \mathbb{R}^d $, where $x' \in \mathbb{R}^{d-2}$ and $x_1$, $x_d \in \mathbb{R}$ . Due to the invariance of $\sigma_{\alpha}$ and partition of unity we just consider a neighborhood of the vectors $\displaystyle x^0 = (0,0',1)$ and $\displaystyle y^0=(\sin(\alpha),0',\cos(\alpha))$. Consider the local coordinates on $\displaystyle S^{d-1}$ around $\displaystyle x^0$ and $\displaystyle y^0$ respectively,
\begin{center}
$\dps x(u)=  \left ( u, 1-\frac{|u|^2}{2} \right ) + O(|u|^3)$, and \\
$\dps y(v)=  \left ( \sin(\alpha)+ v_1, v', \cos(\alpha) - \tan(\alpha)v_1 - \frac{\sec^3(\alpha)}{2} v^2_1 - \frac{\sec(\alpha)}{2}|v'|^2 \right ) + O(|v|^3)$,
\end{center}
where $\displaystyle u=(u_1,u') \in \mathbb{R}^{d-1}$, $\displaystyle v=(v_1,v') \in \mathbb{R}^{d-1}$ and $\displaystyle |u|$, $\displaystyle |v|< \varepsilon$. Note that $\displaystyle x^0 \cdot y^0 = \cos(\alpha)$, but
\begin{align*}
     x(u) \cdot y(v) - \cos(\alpha) & = u_1 \sin(\alpha) + u_1 v_1 + u' \cdot v'  - \tan(\alpha)v_1 \\
      & - \frac{\sec^3(\alpha)}{2} v^2_1  - \frac{\sec(\alpha)}{2}|v'|^2 -\cos(\alpha) \frac{|u|^2}{2}  + O(|u,v|^3) =0
\end{align*}
where the last equality is possible if and only if
\begin{equation*}
  \left ( \tan(\alpha)-u_1 \right )v_1 + \frac{\sec^3(\alpha)}{2} v^2_1 = u_1 \sin(\alpha) + u' \cdot v' - \frac{\sec(\alpha)}{2}|v'|^2 -\cos(\alpha) \frac{|u|^2}{2}  + O(|u,v|^3).
\end{equation*}
By the implicit function theorem in one variable,
\begin{center}
$\displaystyle a_1s+a_2s^2=t \Rightarrow s = a^{-1}_1 t - a^{-3}_1 a_2 t^2 + O(t^3)$ , as $s,t \to 0$.
\end{center}
Then, solving $v_1$ in terms of $u_1$, $u'$ and $v'$ we have
\begin{multline*}
  v_1 = \frac{1}{\tan(\alpha)} \left ( \sum_{k=0}^{\infty} \left ( \frac{u_1}{\tan(\alpha)}\right )^k \right ) \left ( u_1 \sin(\alpha) + u' \cdot v' - \frac{\sec(\alpha)}{2}|v'|^2 -\cos(\alpha) \frac{|u|^2}{2} \right ) \\ -  \frac{1}{\tan^3(\alpha)} \left ( \sum_{k=0}^{\infty} \left ( \frac{u_1}{\tan(\alpha)}\right )^k \right )^3 \left ( \frac{\sec^3(\alpha)}{2} \right ) \left ( u_1 \sin(\alpha) + u' \cdot v' - \frac{\sec(\alpha)}{2}|v'|^2 -\cos(\alpha) \frac{|u|^2}{2} \right )^2
\end{multline*}
\benn
  v_1 = \cos(\alpha)u_1 + \cot(\alpha) (u' \cdot v') - \frac{1}{2\sin(\alpha)} |v'|^2 - \frac{\cos^2(\alpha)}{2\sin(\alpha)} |u'|^2 - \frac{\sin(\alpha)}{2} u^2_1  + O(|u,v|^3).
\eenn
Thus, a neighborhood of $\displaystyle (x^0,y^0)$ can be parameterized as
\begin{align*}
      U_{\alpha} = & \left ( u_1, u', 1-\frac{|u|^2}{2} ; \sin(\alpha)+ \cos(\alpha)u_1 + \cot(\alpha) (u' \cdot v') - \frac{1}{2\sin(\alpha)} |v'|^2 - \frac{\cos^2(\alpha)}{2\sin(\alpha)} |u'|^2  \right. \\
     & \left. - \frac{\sin(\alpha)}{2} u^2_1, v', \cos(\alpha) - \sin(\alpha)u_1 - u' \cdot v' - \frac{\cos(\alpha)}{2} |u'|^2  - \frac{\cos(\alpha)}{2} u^2_1  \right )
\end{align*}
modulo $O(|u,v|^3)$. Then, the Fourier transform of the measure $\sigma_{\alpha}$ can be written as
\benn
    \widehat{\sigma}_{\alpha}(\xi,\eta) = \int e^{-2 \pi i  U_{\alpha} \cdot (\xi, \eta)} \df u_1 \df u'\df v',
\eenn
by invariance assume $(u,u',v')=(0,0',0')$ is a critical point, if there is any. Thus, the gradient of the phase function $\displaystyle U_{\alpha} \cdot (\xi, \eta)$ is given by
\begin{align*}
     \frac{\partial }{\partial u_1} & =  \xi_1 - u_1\xi_d + \cos(\alpha) \eta_1 - \sin(\alpha)u_1 \eta_1 - \sin(\alpha) \eta_d - \cos(\alpha) u_1 \eta_d
\end{align*}
\begin{align*}
     \frac{\partial }{\partial u'} & =  \xi' - u'\xi_d + \cot(\alpha) v'\eta_1 - \frac{\cos^2(\alpha)}{\sin(\alpha)}u' \eta_1 - v' \eta_d + \cos(\alpha) u' \eta_d
\end{align*}
\begin{align*}
     \frac{\partial }{\partial v'} & =  \cot(\alpha)u'\eta_1 - \frac{1}{\sin(\alpha)}v'\eta_1 + \eta' - u' \eta_d
\end{align*}
and the Hessian is given by
\begin{center}
$\dps \left (  - \xi_d - \sin(\alpha) \eta_1 - \cos(\alpha) \eta_d  \right ) \bigoplus^{d-1}_2 \begin{pmatrix} - \xi_d - \frac{\cos^2(\alpha)}{\sin(\alpha)} \eta_1 + \cos(\alpha) \eta_d & \cot(\alpha) \eta_1 - \eta_d \\ \cot(\alpha) \eta_1 - \eta_d & - \frac{1}{\sin(\alpha)} \eta_1 \end{pmatrix}$.
\end{center}
Due to $(u,u',v')=(0,0',0')$ is a critical point, it follows that
\begin{equation}\label{eqn:eq4}
  \xi_1 + \cos(\alpha) \eta_1 - \sin(\alpha) \eta_d =0 \text{ \ and \ }  \xi' = \eta' =0.
\end{equation}
This implies $\dps - \xi_d - \frac{\cos^2(\alpha)}{\sin(\alpha)} \eta_1 + \cos(\alpha) \eta_d = - \xi_d + \frac{\cos(\alpha)}{\sin(\alpha)} \xi_1 $ and $\dps  \cot(\alpha) \eta_1 - \eta_d  =  - \frac{1}{\sin(\alpha)} \xi_1$. Thus the Hessian is given by
\begin{center}
$\dps \left (  - \xi_d - \sin(\alpha) \eta_1 - \cos(\alpha) \eta_d  \right ) \bigoplus^{d-1}_2 \begin{pmatrix}  \xi_d + \frac{\cos(\alpha)}{\sin(\alpha)} \xi_1 $ and $\dps  \cot(\alpha) \eta_1 - \eta_d & \ - \frac{1}{\sin(\alpha)} \xi_1 \\  - \frac{1}{\sin(\alpha)} \xi_1 & - \frac{1}{\sin(\alpha)} \eta_1 \end{pmatrix}$.
\end{center}
Therefore the determinant of the Hessian is given by
\begin{multline*}
         \abs{ - \xi_d - \sin(\alpha) \eta_1 - \cos(\alpha) \eta_d  } \abs{ \pa{ - \xi_d + \frac{\cos(\alpha)}{\sin(\alpha)} \xi_1 } \left ( - \frac{1}{\sin(\alpha)} \eta_1 \right ) - \frac{1}{\sin^2(\alpha)} \xi^2_1 }^{d-2} \\ = \abs{ \xi_d + \sin(\alpha) \eta_1 + \cos(\alpha) \eta_d  } \abs{    \frac{1}{\sin(\alpha)} \eta_1\xi_d - \frac{1}{\sin(\alpha)} \xi_1 \pa{\frac{\cos(\alpha)}{\sin(\alpha)} \eta_1 + \frac{1}{\sin(\alpha)} \xi_1 } }^{d-2}
\end{multline*}
by (\ref{eqn:eq4})
\begin{align*}
      = & \left | \xi_d + \sin(\alpha) \eta_1 + \cos(\alpha) \eta_d  \right | \left |   \frac{1}{\sin(\alpha)} \eta_1\xi_d - \frac{1}{\sin(\alpha)} \xi_1 \eta_d \right |^{d-2}.
\end{align*}
Due to $\xi'=\eta'=0$ we might assume that $\displaystyle \xi,\eta \in \mathbb{R}^2$, thus the last expression is equal to
\begin{align*}
      C \left | \frac{1}{\sin(\alpha)} \right |^{d-2} \left | \xi_d + \sin(\alpha) \eta_1 + \cos(\alpha) \eta_d  \right | |\xi|^{d-2} |\eta|^{d-2} \left | \sin(< \xi, \eta >) \right |^{d-2}.
\end{align*}
Let (\ref{eqn:eq4}) be the first component of a vector in $\mathbb{R}^2$ and let the first factor of the determinant of the Hessian be the second component of such vector,
\begin{center}
  $\displaystyle  \begin{pmatrix} \xi_1 + \cos(\alpha) \eta_1 - \sin(\alpha) \eta_d \\ \xi_d + \sin(\alpha) \eta_1 + \cos(\alpha) \eta_d \end{pmatrix} = \xi + g_{\alpha} \eta$,
\end{center}
where $g_{\alpha} \in O(d)$ is some rotation by $\alpha$. Thus, by stationary phase (see \cite{Wolff2}) we have
\begin{equation*}
  |\h{\sigma}(\xi,\eta)| \leq C_1(\alpha) |\xi + g_{\alpha}\eta|^{-\frac{1}{2}} |\xi|^{-\frac{(d-2)}{2}} |\eta|^{-\frac{(d-2)}{2}} | \sin( < \xi , \eta >)| ^{-\frac{(d-2)}{2}}.
\end{equation*}
Where $\dps C_1(\alpha)= C  \abs{ \frac{1}{\sin(\alpha)} }^{-\frac{(d-2)}{2}}$.
\vspace{2mm}


\section{Proof of Lemma \ref{lem2.3}}

Case 1: If $1 \lesssim <\xi,\eta>$ and $R2^j \approx |\xi| \lesssim |\xi +g_{\alpha}\eta|$ then the result is immediate.\\

Case 2: If $ <\xi,\eta> $ is small and $(R2^j)^{-1} \approx |\xi|^{-1} \approx |\xi +g_{\alpha}\eta|^{-1}$. By using polar coordinates we have
\begin{multline*}
   \int \limits_{  \begin{matrix} |\xi| \approx R2^j \\ < \xi,\eta > \text{\ small} \end{matrix} } |\xi + g_{\alpha}\eta|^{-1} (\sin( < \xi , \eta > )) ^{-(d-2)} \df \xi  \\  \lesssim  \pa{R2^j }^{-1+d} \int \limits_{S^{d-1} } (\sin( < \omega , \eta > )) ^{-(d-2)} \df \omega_{d-1} = 2 \pi (R2^j)^{-1+d}.
\end{multline*}

 Case 3: $ |\xi +g_{\alpha}\eta| $ is small.\\
Due to $|\xi| \approx |\eta|$, then $1 \lesssim <\xi,\eta>$. Thus
\begin{multline*}
    \int \limits_{  \begin{matrix} |\xi| \approx R2^j \\ |\xi + g_{\alpha}\eta| \text{\ small} \end{matrix} } |\xi + g_{\alpha}\eta|^{-1} (\sin( < \xi , \eta > )) ^{-(d-2)} \df \xi  \lesssim \int \limits_{  |\zeta| \lesssim R2^j } |\zeta|^{-1} \df \zeta \lesssim (R2^j)^{-1+d},
\end{multline*}
where $\zeta = \xi + g_{\alpha}\eta$.
\vspace{2mm}

\section{Proof of Lemma \ref{PHprinc}}

A result similar to Lemma \ref{PHprinc} was given by Iosevich, Mourgoglou and Senger (See Lemma $2.1$ in \cite{IosevichMourgoglouSenger}) , where the authors show the existence of two subsets that satisfy $(i)$ and $(ii)$. We will adapt their proof and see that $(iii)$ follows immediately from the construction (see Figure \ref{fig:M1}) given below.\\

We will use a stopping time argument. Consider the unit cube $\left [ 0,1 \right ]^d$ in $\mathbb{R}^d$, and let $C_{\mu}$ be the constant in the Frostman condition $\displaystyle \mu(B(x,r)) \leq C_{\mu} r^s$, for all $x \in \mathbb{R}^d$ and $r>0$. Moreover, without lost of generality we may assume that $C_{\mu}>1$.\\

Lets subdivide the unit cube $\left [ 0,1 \right ]^d$ into $4^d$ smaller cubes with side-length $1/4$. Let $\Omega_k$, $1 \leq k \leq 2^d$, be the collection of $2^d$ sub-cubes such that no two cubes of the same collection touch each other. By pigeon hole principle at least one of the collections $\Omega_k$ has measure greater or equal to $ \frac{1}{2^d}$, that is,
\begin{center}
  $\displaystyle \mu \left ( \bigcup_{Q \in \Omega_k} Q  \right ) \geq \frac{1}{2^d}$ for some $k$.
\end{center}
Thus we have the following cases
\begin{itemize}
  \item [(1)] If there are three cubes in the collection $\Omega_k$, say $Q_1$, $Q'_1$ and $Q''_1$, such that $\mu (Q_1)$, $\mu (Q'_1)$ and $\mu (Q''_1)$ are greater or equal to $\frac{c}{2^d}$, for some $c>0$ , then we are done.
\end{itemize}

\begin{itemize}
  \item [(2)] If there are only two cubes in the collection $\Omega_k$, say $Q_1$ and $Q'_1$ such that both $\mu (Q_1)$ and $\mu (Q'_1) $ are greater or equal to $\frac{c}{2^d}$, for some $c>0$, then we have the following subcases
      \begin{itemize}
        \item [(a)] If $c \geq 1$, then clearly we have $\mu (Q_1) \geq \frac{1}{2^d} $, and we repeat the same procedure on the cube $Q_1$.
        \item [(b)] If $c < 1 $, then for at least one of the cubes, say $Q_1$ we must have $\mu (Q_1) \geq \frac{1}{2^{d+1}} $, and we repeat the same procedure on the cube $Q_1$. Otherwise, if we had $0 < \mu (Q_1) , \mu (Q'_1) < \frac{1}{2^{d+1}}$, then
            \begin{center}
              $\displaystyle \mu \left ( \bigcup_{Q \in \Omega_k \setminus \left \{ Q_1,Q'_1 \right \} } Q  \right ) >0 $.
            \end{center}
            Thus, there must be a cube in $\Omega_k \setminus \left \{ Q_1,Q'_1 \right \}$ with positive $\mu$ measure, and we have case $(1)$.
      \end{itemize}
\end{itemize}

\begin{itemize}
  \item [(3)] If there is only one cube in the collection $\Omega_k$, say $Q_1$, with positive measure then we must have $\mu (Q_1) \geq \frac{1}{2^d}$. Then, we repeat the procedure on the cube $Q_1$.
\end{itemize}
 Note that we can repeat this process, and at each stage check for three cubes. Lemma $2.1$ in \cite{IosevichMourgoglouSenger} showed that we cannot fail to find two cubes with positive measure, for the sake of completeness we reproduce their proof in claim 1. \\

 {\bf Claim 1:} There is a collection in which we can find two cubes with positive measure.\\
 Suppose that at every iteration we cannot find two cubes with positive measure. Then at the $n-$th iteration, we have a cube $Q_n$ of side length $\frac{1}{4^{n}}$ for which $\mu(Q_n) \geq \frac{1}{2^{dn}}$. By the Frostman measure condition we have $\dps \frac{1}{2^{dn}} \leq \mu(Q_n) \leq C_{\mu} \frac{1}{4^{ns}}$, from which we obtain $\dps n \leq \frac{\log_2(C_{\mu})}{2s-d}$ for every $n$ which is a contradiction.

 {\bf Claim 2:} There is a collection in which we can find three cubes with positive measure. Suppose that at every iteration we cannot find three cubes with positive measure. Because of Claim $1$, we just study case $(2b)$. If we fail to find a third cube at the $n-$th iteration, we obtain a cube $Q_n$ of side-length $\frac{1}{4^{n}}$ for which $\mu(Q_n) \geq \frac{1}{2^{(d+1)n+1}}$. By the Frostman measure condition we have $\dps \frac{1}{2^{(d+1)n+1}} \leq \mu(Q_n) \leq C_{\mu} \frac{1}{4^{ns}}$, from which we obtain $\dps n \leq \frac{\log_2(C_{\mu}) +1}{2s-(d+1)}$ for every $n$ which is a contradiction.\\

\begin{center}
  \begin{tikzpicture}[scale=0.6]
\draw[step=1cm,white,very thin] (0,0) grid (16,16);
\draw[step=1cm,gray, thin] (8,8) grid (12,12);
\fill[blue!40!white] (8,9) rectangle (9,10);
\fill[blue!40!white] (10,9) rectangle (11,10);
\fill[blue!40!white] (10,11) rectangle (11,12);
\draw[step=4cm,black, thick ] (0,0) grid (16,16);
\draw (0,0) rectangle (16,16);
\draw (2,14) node { \verb!1!};
\draw (2,6) node {  \verb!1!};
\draw (10,14) node { \verb!1!};
\draw (10,6) node { \verb!1!};
\draw (6,14) node {\verb!2!};
\draw (14,14) node {\verb!2!};
\draw (14,6) node {\verb!2!};
\draw (6,6) node {\verb!2!};
\draw (2,2) node {\verb!3!};
\draw (2,10) node {\verb!3!};
\draw (10,2) node {\verb!3!};
\draw (6,2) node {\verb!4!};
\draw (14,2) node {\verb!4!};
\draw (14,10) node {\verb!4!};
\draw (6,10) node {\verb!4!};
\draw (10.5,11.5) node {$E_1$};
\draw (8.5,9.5) node {$E_3$};
\draw (10.5,9.5) node {$E_2$};
\end{tikzpicture}
\captionof{figure}{When $d=2$, this picture shows the case in which it is possible to find three cubes with positive $\mu-$measure in the second iteration. In the first iteration we have four collections of four cubes, which was not enough to obtained three cubes with positive $\mu-$measure, but the upper-right cube in the third collection has positive $\mu-$measure. By subdividing this cube into four collections of four smaller cubes we succeed to find three cubes with positive $\mu-$measure. This picture is based on the graphical description given in \cite{IosevichMourgoglouSenger}.}\label{fig:M1}
\end{center}

Note that from Figure \ref{fig:M1} we conclude that it is possible to find sets $E_1$, $E_2$ and $E_3$ with positive $\mu-$measure such that $\dps 0 < c_3 \leq \alpha(x,z,y) \leq c_4 < \frac{\pi}{2} $, where $x \in E_1$, $y \in E_2$ and $z \in E_3$. The reader might be concerned, for instance, in the case where $\mu$ is concentrated in the bottom left corner of $E_1$, top right corner of $E_3$ and bottom left corner of $E_2$. In that case we can subdivide each of the cubes $E_1$, $E_2$ and $E_3$ into $4^d$ smaller cubes, and considering the bottom left cube that is contained in $E_1$, the top right cube that is contained in $E_3$ and the bottom left cube that is contained in $E_2$ (See Figure \ref{fig:M2} below). Thus we have sets $E_{1,1}$, $E_{2,1}$ and $E_{3,1}$ such that $\dps 0 < c_3 \leq \alpha(x,z,y) \leq c_4 < \frac{\pi}{2} $, where $x \in E_{3,1}$, $y \in E_{2,1}$ and $z \in E_{1,1}$.
\begin{center}
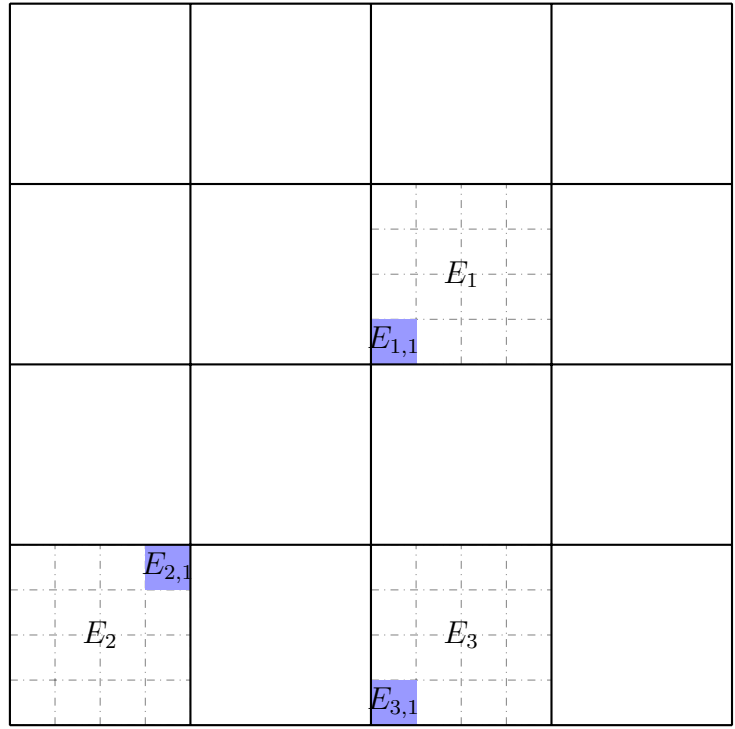

\begin{tikzpicture}[scale=0.6]
\draw[step=1cm,white,very thin] (0,0) grid (16,16);
\draw[step=1cm,dashdotted,gray,very thin] (8,8) grid (12,12);
\draw[step=1cm,dashdotted,gray,very thin] (0,0) grid (4,4);
\draw[step=1cm,dashdotted,gray,very thin] (8,0) grid (12,4);
\fill[blue!40!white] (8,8) rectangle (9,9);
\fill[blue!40!white] (3,3) rectangle (4,4);
\fill[blue!40!white] (8,0) rectangle (9,1);
\draw[step=4cm,black, thick ] (0,0) grid (16,16);
\draw (0,0) rectangle (16,16);
\draw (2,14) node { };
\draw (2,6) node {  };
\draw (10,14) node { };
\draw (10,6) node { };
\draw (6,14) node {};
\draw (14,14) node {};
\draw (14,6) node {};
\draw (6,6) node {};
\draw (2,2) node {};
\draw (2,10) node {};
\draw (10,2) node {};
\draw (6,2) node {};
\draw (14,2) node {};
\draw (14,10) node {};
\draw (6,10) node {};
\draw (10,10) node {$E_{1}$};
\draw (2,2) node {$E_{2}$};
\draw (10,2) node {$E_{3}$};
\draw (8.5,8.5) node {$E_{1,1}$};
\draw (3.5,3.5) node {$E_{2,1}$};
\draw (8.5,0.5) node {$E_{3,1}$};
\end{tikzpicture}
\captionof{figure}{By taking $x \in E_{3,1}$, $y \in E_{2,1}$ and $z \in E_{1,1}$ we have $\dps 0 < c_3 \leq \alpha(x,z,y) \leq c_4 < \frac{\pi}{2} $}\label{fig:M2}
\end{center}
\vspace{2mm}
Any other case in which we might have the same issue can be fixed by considering smaller cubes or by changing the order in which the cubes are being considered.


\section{Proof of Corollary \ref{incthm} }

Corollary \ref{incthm} follows immediately from the equation
    \begin{multline*}
    \mu \times \mu \times \mu  \left ( \left \{ (x,y,z): t \leq |x-z| \leq t + \ve ,\ r \leq |y-z| \leq r + \ve, \right. \right. \\
    \left. \left.  \alpha \leq \alpha(x,y,z) \leq \alpha + \ve  \right \} \right ) = \int_{t}^{t + \ve} \int_{r}^{r + \ve} \int_{\alpha}^{\alpha + \ve} \delta(\mu) (\widetilde{t},\widetilde{r},\widetilde{\alpha}) \df \widetilde{t} \df \widetilde{r} \df \widetilde{\alpha},
    \end{multline*}
and the fact that the density of the measure $\df \delta(\mu)$ is continuous and compactly supported.


\bigskip

\end{document}